\documentclass [12pt] {article}
\usepackage{amsfonts}
\usepackage{amssymb}

\newcommand{\qed} {\hspace {0.1in} \rule {1.5mm} {3.5mm}}

\newtheorem{theorem}{Theorem}

\newtheorem{proposition}{Proposition}[section]
\newtheorem{definition}{Definition}[section]

\def\ek{\mbox{End}_k}

\def\dim{{\rm dim}}

\def\<{\langle}
\def\>{\rangle}

\def\proof{\smallskip\noindent{\bf Proof:} }

\def\ca{\mbox{$\cal A$}}
\def\cb{\mbox{$\cal B$}}

\def\Ran{\,\mbox{Ran}\,}
\def\Ker{\,\mbox{Ker}\,}
\def\Mat{\,\mbox{Mat}\,}

\def\to{\rightarrow}

\title{On Algebras That Almost Have Finite Dimensional Representations} 
\author{{\sc G\'abor Elek}
\cr Mathematical Institute of
the Hungarian Academy of Sciences\cr P.O. Box 127, H-1364 Budapest, Hungary\cr
elek@renyi.hu}
\date{}
\begin{document}

\maketitle
\noindent{\bf Abstract.} We introduce the notion of almost finite
dimensional representability of algebras and study its connection
 with the classical finiteness
conditions.
\vskip 0.2in
\newpage
\section{Introduction}
In \cite{VG}  Gordon and Vershik introduced the notion of $LEF$-groups. A group
$G$ is called $LEF$ (locally embeddable into finite groups) if for any finite
subset $1\in S\subset G$ there exists an injective map
$\psi_S:S\to G_S$ into a finite group $G_S$ such that:
$\psi_S(1)=1$ and $\psi_S(ab)=\psi_S(a)\cdot\psi_S(b)$, whenever
$a, b, ab\in S$.

\noindent
In the same paper they defined $LEF$-algebras as well. Let $k$ be a
commutative
field and $\ca$ be $k$-algebra. Then $\ca$ is called $LEF$ if for any finite
dimensional subspace $L\subset \ca$ there exists an injective linear map
$\psi_L:L\to \ek(V)$ for some finite dimensional $k$-space $V$, such that:
$\psi_L(1)=1$ and $\psi_L(xy)=\psi_L(x)\psi_L(y)$, whenever $x,y,xy\in L$.

\noindent
Let us recall the three class of pathological rings from Cohn's classical
paper \cite{Cohn}.
\begin{description}
\item[I] There exist $A\in \Mat_{m\times n}(R)$ and
 $B\in \Mat_{n\times m}(R)$, $m\neq n$,
 such that $AB=I$ and $BA=I$  (in other words, the ring $R$ does not have the 
 {\it IBN property}).

\item[II] There exist $A\in \Mat_{m\times n}(R)$, and
 $B\in \Mat_{n\times m}(R)$, $m>n$, such that
 $AB=I$ (or
  the ring $R$ does not satisfy the {\it rank condition}).

\item[III] There exist $A\in \Mat_{n\times n}(R)$, and
 $B\in \Mat_{n\times n}(R)$, such that $AB=I$ and $BA\neq I$ (or
 the ring $R$ is not {\it stably finite}).
\end{description}

\noindent
Obviously, stable finiteness implies the rank condition, and the rank
condition in turn implies the IBN property. Cohn  showed that this categories
are strictly containing each other. Now, let us observe that $LEF$-algebras
can not be pathological at all.
\begin{proposition}
$LEF$-algebras are stably finite.
\end{proposition}
\proof Let $\ca$ be a $LEF$-algebra.
Suppose that 
$A\in \Mat_{n\times n}(\ca)$, 
 $B\in \Mat_{n\times n}(\ca)$, such that $AB=I$. Consider the finite
dimensional subspace $L$ generated by $1\in\ca$, the entries of $A$ and $B$ and
all the possible products in the form of $xy$, where $x$ is an entry of
$A$ and $y$ is an entry of $B$. Let us consider an injective map
$\psi_L:L\to \ek(V)$ as in the definition of $LEF$-algebras. Obviously
it extends to a map $\widehat{\psi}_L:\Mat_{n\times n}( L)\to \ek(V^n)$
such that $\widehat{\psi}_L(A)\widehat{\psi}_L(B)=I\in \ek(V^n)$ and
$\widehat{\psi}_L(B)\widehat{\psi}_L(A)\neq I$, where 
$\Mat_{n\times n}(L)$ is the finite dimensional vector space of matrices
with entries from $L$. The contradiction is clear, since the finite
dimensional
matrix algebras are stably finite. \qed

\noindent
There is a rather simple ring theoretical obstruction as well
that excludes, say, all the Weyl-algebras from the category of
$LEF$-algebras.
\begin{proposition}
Let $\ca$ be a finitely presented infinite dimensional {\it simple}
affine algebra. Then $\ca$ is not a $LEF$-algebra.
\end{proposition}
\proof
Let $a_1,a_2,\dots, a_n$ be generators of $\ca$ and let
$f_i(a_1,a_2,\dots,a_n)=0$ be presentations, $1\leq i\leq l$, where
the $f_i$'s are polynomials of non-commuting variables.
Easy to see that there exists a finite dimensional subspace $L\subset \ca$ and
an injective linear map $\psi_L:L\to \ek(V)$
into a finite dimensional matrix algebra over $k$ such that;
$$f_i(\psi_L(a_1),\psi_L(a_2),\dots,\psi_L(a_n))=0,\quad 1\leq i\leq l\,.$$
Therefore the algebra generated by
the matrices $\psi_L(a_1),\psi_L(a_2),\dots \psi_L(a_n)$ is
a non-trivial quotient of $\ca$, this leads to a contradiction. \qed

\noindent
In \cite{Gro}, Gromov defined a weaker version of the $LEF$-property for
groups, that he called initial amenability. 
This led us to the following definition.
\begin{definition}\label{alm}
A unital $k$-algebra $\ca$ almost has finite dimensional representations if
for any finite dimensional subspace $1\in L\subset \ca$ and $\epsilon>0$,
there exists a finite dimensional vector space $V$ together with a subspace
$V_\epsilon\subset V$ such that 
\begin{itemize}
\item There exists a linear (not necessarily injective) map
$\psi_{L,\epsilon}:L\to \ek(V)$ such that $\psi_{L,\epsilon}(1)=I$ and
$\psi_{L,\epsilon}(a)\psi_{L,\epsilon}(b)( v)= \psi_{L,\epsilon}(ab) (v)$
whenever $a,b,ab\in L$ and $v\in V_\epsilon$.
\item $\frac {\dim_k V-\dim_k V_\epsilon}{\dim_kV} <\epsilon\,.$
\end{itemize}
\end{definition}
We shall call such maps $\epsilon$-almost representations of $L$.
Obviously algebras (e.g. group algebras) that actually have  finite
dimensional representations or  $LEF$-algebras are in the class above.
Clearly, if $\ca$ almost has finite dimensional representations so does any
subalgebra $\cb\subset\ca$ containing the unit of $\ca$.
Also, if $\ca$ almost has finite dimensional representations and $\ca$ is a 
quotient of the algebra $\cb$, then $\cb$ almost has finite dimensional
representations. 
Note that we shall give examples of simple infinite dimensional algebras
almost having finite dimensional representations.
Our main results are the following two theorems.
\begin{theorem}

\noindent
\begin{enumerate} \item
Let $\ca$ and $\cb$ be $k$-algebras almost having finite dimensional
representation, then both $\Mat_{n\times n} \ca$ and $A\otimes_k \cb$ are algebras
almost having finite dimensional representations.
\item
If $\ca$ almost has finite dimensional representations then it satisfies
the rank condition. 
\end{enumerate}\end{theorem}
\begin{theorem}
If $\ca$ is a simple algebra almost having finite dimensional representations
then it is stably finite.
\end{theorem}
Actually, we shall construct an ideal $RR(\ca)$ in any
 $k$-algebra $\ca$
such that $RR(\ca)\neq\ca$ if and only if $\ca$ almost has a finite
dimensional representation and if $RR(\ca)=0$ then $\ca$ is stably finite.
For any pair of algebras $\ca$,$\cb$ and any algebra homomorphism
 $\tau:\ca\to\cb$, we shall see that $\tau(RR(\ca))\subset RR(\cb)$.
\section{Examples}
In \cite{Row}, Rowen proved that any algebra of subexponential growth
is an IBN-algebra. In \cite{Elek}, we extended Rowen's result to the
larger class of {\it amenable algebras}.
\begin{definition}
The $k$-algebra $\ca$ is amenable if for any finite dimensional linear
subspace $1\in B\subset \ca$ and $\epsilon>0$ there exists
a non-trivial finite dimensional linear subspace $Q\subset A$ such that
$$\frac {\dim_k BQ-\dim_k Q}{\dim_k Q} <\epsilon\,.$$
Here $BQ$ denotes the subspace spanned by the vectors in the form $bq$,
where $b\in B$ and $q\in Q$.
\end{definition}
\begin{proposition}\label{ame}
If $\ca$ is amenable, then it almost has finite dimensional representations.
\end{proposition}
\proof
Fix a finite dimensional subspace $1\in L\subset \ca$. Let $\{Q_n\}_{n\geq 1}$
be a sequence of finite dimensional subspaces in $\ca$ such that
$$\frac {\dim_k LQ_n-\dim_k Q_n}{\dim_k Q_n} <\frac {1} {n}\, $$
and consider linear subspaces $T_n$ such that $Q_n\oplus T_n=\ca$ as vector
spaces. Denote the projection onto $Q_n$ by $P_n$.
Now, let $\psi_n:L\to \ek(Q_n)$ be defined as
$$\psi_n(a) v=P_n(av)\,.$$
Note that $\psi_n(ab)(v)=\psi_n(a)\psi_n(b)(v)$ if
$P_n(bv)=bv$, that is $v\in \Ker(m_b-P_nm_b)\,$
where $m_b$ is the left-multiplication by the element $b$.
Observe that $\dim_k \Ran(m_b-P_nm_b)\leq \frac{1}{n}\dim_k Q_n\,.$
Hence
$$\dim_k \Ker(m_b-P_nm_b)\geq \frac{n-1}{n} \dim_k Q_n\,.$$
Consequently, for any $0<\epsilon<1$,
$$\dim_k\bigcap_{x\in L} \Ker(m_x-P_n m_x)\geq (1-\epsilon)
\dim_k Q_n\,,$$
if $n$ is large enough. Let $Q^\epsilon_n=\bigcap_{x\in L} \Ker(m_x-P_n m_x)$.
Then, if $v\in Q^\epsilon_n$:
$$\psi_n(xy)(v)=\psi_n(x)\psi_n(y)(v)\,,$$
provided that $x,y,xy\in L$. Also,
$$\frac {\dim_k Q_n-\dim_k Q^\epsilon_n}{\dim_k Q_n} <\epsilon\,. $$
Hence $\ca$ almost has finite dimensional representations. \qed

\noindent
Note that by Theorem 1. we can immediately see that amenable algebras
satisfy the rank condition. 

\noindent
There is however an even more direct relation between amenability and
the property studied in our paper. Let $G(V,E)$ be an infinite
connected graph of bounded vertex degree. Then $V(G)$ is a discrete metric 
space with the shortest path-metric $d_G$. The graph is called {\it amenable}
if its isoperimetric constant is zero, or in other words, there exists
finite subsets $F_n\subset V(G)$ such that for any $k\geq 1$:
$$\lim_{n\to\infty}\frac{|B_k(F_n)|}{|F_n|}=1\,,$$ where $B_k(F_n)$
is the $k$-neighbourhood of $F_n$ in the metric $d_G$:
$$B_k(F_n)=\{y\in V :\,\mbox{there exists}\,\,x\in F_n;\,d_G(x,y)\leq k\}\,.$$
The translation algebra $T_k(G)$ of $G$ is the set of all square matrices
$A$ indexed by $V\times V$ with entries from $k$ such
that $A(x,y)=0$, whenever $d(x,y)>A_l$ for some constant depending only $A$
\cite{AMP}.
The following result is motivated by Proposition 4.1 \cite{AMP}.
\begin{proposition}
The graph $G$ is amenable if and only if the translation algebra $T_k(G)$
almost has finite dimensional representations.
\end{proposition}
\proof
Suppose that $G$ is amenable. Let $Q_n$ be the vectorspace over $k$ spanned
by the elements of $F_n$ as a formal basis. Then for any finite dimensional
subspace $1\in L\subset T_k(G)$, clearly
$$\lim_{n\to 0} 
\frac {\dim_k LQ_n-\dim_k Q_n}{\dim_k Q_n}=0\,.$$
Thus one can use exactly the same argument as in Proposition \ref{ame}
to show that $T_k(G)$ almost has finite dimensional representations.
If $G$ is non-amenable, then by \cite{DSS} (see also \cite{Elek2})
there exists a partition $V=V_1\cup V_2$ and bijective maps 
$\phi_1:V\to V_1$, $\phi_2:V\to V_2$ such that
\begin{equation}
\label{e2}
\sup_{x\in V} d_G(x,\phi_1(x))<\infty, \quad
\sup_{x\in V} d_G(x,\phi_2(x))<\infty
\end{equation}
Let $A(x,y)=1$ if $x=\psi_1(y)$, otherwise let $A(x,y)=0$. Similarly,
let $B(x,y)=1$ if $x=\psi_2(y)$, otherwise let $B(x,y)=0$.
By (\ref{e2}) the matrices $A$ and $B$ are in $T_k(G)$.
Obviously, $AA^T=I$, $BB^T=I$ and $A^TA+B^TB=I$. Therefore the algebra
generated by $A,A^T,B,B^T$ is not an IBN-algebra, hence $T_k(G)$ is
not an IBN-algebra as well. By Theorem 1. the proposition follows. \qed

\section{The proof of Theorem 1.}
Suppose that $\ca$ almost has finite dimensional representations.
Let $L\subset \Mat_{n\times n}(\ca)$ be a finite
dimensional $k$-space. Then there exists a finite dimensional $k$-space
$M\subset \ca$ such that $L\subset \Mat_{n\times n} (M)$. 
Let $\psi_{M,\epsilon}:M\to \ek(V)$ be an $\epsilon$-almost representation
of $M$. Then we can extend it to $\widehat{\psi}_{M,\epsilon}:
\Mat_{n\times n}(M)\to \ek(V^n)$ in a natural manner. If $s\in (V_\epsilon)^n$,
then $\widehat{\psi}_{M,\epsilon}(ab) (s)= \widehat{\psi}_{M,\epsilon}(a)
\widehat{\psi}_{M,\epsilon}(b) (s)$. Obviously,
$\lim_{\epsilon\to 0}\frac
{\dim_k V^n-\dim_k (V_\epsilon)^n}{\dim_k V^n}=0\,.$
Hence $\Mat_{n\times n}(\ca)$ almost has finite dimensional representations.

\noindent
Now let $\ca$ and $\cb$ be algebras almost having finite dimensional representations.
If $L\subset \ca\otimes_k \cb$ is a finite dimensional vectorspace, then
$L\subset M\otimes_k N$, where $M\subset \ca$ and $N\subset \cb$ are finite
dimensional vectorspaces. If we choose $\epsilon$-almost representations
$\psi_{M,\epsilon}:M\to \ek(V)$, $\psi_{N,\epsilon}:N\to \ek(W)$
, then we obtain a linear map
$$\psi_{M,\epsilon}\otimes \psi_{N,\epsilon}:M\otimes N\to \ek(V)\otimes
\ek(W)\sim \ek(V\oplus W)\,.$$
Obviously,
$$(\psi_{M,\epsilon}\otimes \psi_{N,\epsilon})(ac\otimes bd)(v,w)=
(\psi_{M,\epsilon}\otimes \psi_{N,\epsilon})(a\otimes b)\cdot(
\psi_{M,\epsilon}\otimes \psi_{N,\epsilon})(c\otimes d) (v,w)\,,$$
provided that $v\in V_\epsilon$ and $w\in W_\epsilon$.
Consequently, $\ca\otimes_k \cb$ almost has finite dimensional representations.

\noindent
Finally, suppose that $A\in \Mat_{m\times n}(\ca)$, $B\in \Mat_{n\times
m}(\ca)$, $m>n$
and $AB=I$. Then, for any $\epsilon$ we have finite dimensional
$k$-spaces $V$, $V_\epsilon,$ 
$\frac {\dim_k V-\dim_k V_\epsilon}{\dim_kV} <\epsilon$ and linear maps
$A_\epsilon:V^n\to V^m$, $B_\epsilon:V^m\to V^n$ such that
$A_\epsilon B_\epsilon (v)=v$, if $v\in V^m_\epsilon$. Note however
that the range of $A_\epsilon B_\epsilon$ is at most $n\cdot \dim_k V$-
dimensional. However, if $\epsilon$ is small enough,
$$m\cdot \dim_k V_\epsilon> n\cdot \dim_k V\,\,,$$
leading to a contradiction. \qed
\section{The rank radical}
In order to prove Theorem 2. we need the notion of a rank radical of an
algebra, $RR(\ca)$.
\begin{definition}
If $\ca$ is not an algebra that almost has a finite dimensional
representation, set $RR(\ca)=\ca$. Otherwise, let $p\in RR(\ca)$
if there exists a finite dimensional subspace $\{1,p\}\subset L\subset \ca$
such that for any $\delta>0$ there exists $n_\delta>0$ with the
following property:

\noindent
If $0<\epsilon < n_\delta$ and $\psi_{L,\epsilon}:L\to \ek(V)$
is an $\epsilon$-almost representation
then
$$\dim_k \Ran\,(\psi_{L,\epsilon} (p))< \delta\, \dim_k V\,\quad.$$
Note that if $\phi:\ca\to \ek(V)$ is an algebraic homomorphism then
$\phi(RR(\ca))=0$.
\end{definition}
\begin{proposition}\label{p41}
$RR(\ca)$ is an ideal.
\end{proposition}
Let $a\in \ca$, $p\in RR(\ca)$ be arbitrary elements and $\{a,p\}\subset
L \subset \ca$ be the finite dimensional vectorspace of the previous
definition. Let $L\subset L'\subset \ca$ be any finite dimensional
vectorspace that contains $a$ and $ap$.
If $\epsilon< n_{\frac{1}{2} \delta}$:
$$\dim_k \Ran (\psi_{L',\epsilon}(p))\leq \frac{1}{2}\delta\cdot \dim_k V$$
since $L\subset L'$. Note that $\psi_{L',\epsilon}(ap)( v)= 
\psi_{L',\epsilon}(a)\cdot \psi_{L',\epsilon}(p)( v)$
if $v\in V_\epsilon$, thus
$$\dim_k \Ran (\psi_{L',\epsilon}(ap))\leq \frac{1}{2}\delta\cdot \dim_k
+\epsilon\cdot \dim_k V\,\quad.$$
Therefore if $\epsilon$ is small enough then
$$\dim_k \Ran (\psi_{L',\epsilon}(ap))\leq \delta\cdot \dim_k V\,.$$
Hence $ap\in RR(\ca)$. Similarly, $pa\in RR(\ca)$. \qed

\noindent
Observe that if $p\notin RR(\ca)$, then for any finite dimensional
subspace $\{1,p\}\subset L\subset \ca$ and $\epsilon>0$
there exists an $\epsilon$-almost representation such that
$$\dim_k \Ran (\psi_{L,\epsilon}(p))\geq C_L \cdot \dim_k V\,,$$
where $C_L$ depends only on $L$. Hence the following proposition easily
follows :
\begin{proposition}
Let $\ca,\cb$ be $k$-algebras and let
$\tau:\ca\to \cb$ be an algebra homomorphism.
Then $\tau(RR(\ca))\subset RR(\cb)$.
\end{proposition}
Now we prove a statement that immediately implies Theorem 2. by
Proposition \ref{p41}.
\begin{proposition}
Suppose that $RR(\ca)=0$, then $\ca$ is stably finite.
\end{proposition}
\proof
Let $T,S\in \ek(k^l)$, suppose that $V\subseteq k^l$, $\dim_k(V)\geq (1-\epsilon)l$
and $TS(v)=v$ if $v\in V$. Then 
\begin{equation}\label{euj3}
\dim_k \Ran (TS-ST)\leq 2\epsilon l. \end{equation}
Indeed, if $w\in S(V)$, then $ST(w)=w$. Hence if $w\in V\cap S(V)$, then
$(TS-ST)v=0$. Obviously, $\dim_k(V\cap S(V))\geq (1-2\epsilon)l$. 
Thus (\ref{euj3}) follows.
Now suppose that $A\in \Mat_{n\times n} (\ca)$, $B\in \Mat_{n\times n}(\ca),$
 $AB=I$, $BA=P\neq I$. 
Let $L$ be a finite dimensional vectorspace spanned by the elements
in the form $1,x,y,xy,yx$, where $x$ is an entry of $A$ and $y$ is
an entry of $B$. Let $p\in L$ be a non-zero entry of $AB-BA$.

\noindent
Let $\{\psi^m_L:L\to \ek(V_m)\}_{m\geq 1}$ be a sequence of
$\frac{1}{m}$-almost representations such that
\begin{equation} \label{euj4}
\dim_k \Ran (\psi^m_L(p))>C\cdot\dim_k V_m\,
\end{equation}
for some constant $C>0$.
We have the extension:
$$\widehat{\psi}^m_L: \Mat_{n\times n} (L)\to \ek((V_m)^n)$$
such that $\widehat{\psi}^m_L(A)\cdot \widehat{\psi}^m_L(B) (v)=v$
if $v\in W_m\subset (V_m)^n$, where
$$\lim_{m\to\infty} \frac{ n\cdot\dim_k V_m-\dim_k W_m}{n\cdot \dim_k
V_m}=0\,.$$
By (\ref{euj4}) $$ \frac{\dim_k \Ran(\widehat{\psi}^m_L (AB-BA))}
{n\cdot\dim_k V_m}$$ can not converge to zero.
On the other hand, by the observation at the beginning of the proof
$$\lim_{m\to\infty}\frac{\dim_k \Ran ( \widehat{\psi}^m_L (A) \widehat{\psi}^m_L (B)-
\widehat{\psi}^m_L (B) \widehat{\psi}^m_L (A))}
{n\cdot\dim_k V_m}=0$$
in contradiction with the fact that 
$$\widehat{\psi}^m_L (AB-BA) (v)
=(\widehat{\psi}^m_L (A) \widehat{\psi}^m_L (B)-
\widehat{\psi}^m_L (B) \widehat{\psi}^m_L (A)) (v)$$
if $v\in W_m$. \qed

\end{document}